\documentclass[10pt,leqno]{article}
\usepackage{amssymb,amsbsy,amsmath,amsfonts,amssymb,amscd}

\usepackage[english]{babel}
\usepackage[T1]{fontenc}
\usepackage{indentfirst}

\newtheorem{statement}{}
\newtheorem{theoreme}[statement]{Theorem}
\newtheorem{lemme}[statement]{Lemma}

\newtheorem{proposition}[statement]{Proposition}

\newtheorem{corollaire}[statement]{Corollary}

\newcommand\C{\mathbb C}
\newcommand\N{\mathbb N}

\newcommand\T{\mathbb T}
\newcommand\D{\mathbb D}
\newcommand\Z{\mathbb Z}
\newcommand\e{{\rm e}}
\renewcommand\P{\mathbb P}
\newcommand\eps{\varepsilon}
\newcommand\ind{{\rm 1\kern-.30em I}}

\title{Weak compactness and Orlicz spaces}
\author{Pascal Lef\`evre, Daniel Li,\\ Herv\'e Queff\'elec, Luis Rodr{\'\i}guez-Piazza}

\date{\footnotesize \today}

\begin{document}

\maketitle

\noindent{\bf Abstract.} \emph{We give new proofs that some Banach spaces have 
Pe{\l}czy\'nski's property $(V)$.}
\smallskip


\section{Introduction.}

Recall that a Banach space $X$ is said to have Pe{\l}czy\'nski's property $(V)$ if one has a 
good weak-compactness criterion in the dual space $X^\ast$ of $X$, namely: every subset $A$ of 
$X^\ast$ is relatively weakly compact whenever it has the following property (easily seen 
necessary):
\begin{displaymath}
\lim_{n\to +\infty} \sup_{x^\ast\in A} | x^\ast(x_n)| = 0
\end{displaymath}
for every weakly unconditionaly Cauchy sequence $(x_n)_n$ in $X$ ({\it i.e.} such that 
$\sum_{n\geq 1} |x^\ast(x_n)| <\infty$ for any $x^\ast\in X^\ast$). Equivalently, $X$ has 
Pe{\l}czy\'nski's property $(V)$ if and only if for every Banach space $Z$ and every non-weakly 
compact operator $T\colon X\to Z$, there exists a subspace $X_0$, isomorphic to $c_0$, 
such that $T$ is an isomorphism between $X_0$ and $T(X_0)$. Beside the reflexive spaces (and 
in particular the $L^p$ spaces for $1<p<\infty$), the spaces ${\cal C}(S)$ of continuous 
functions on compact spaces $S$ have property $(V)$; in particular $L^\infty$ has $(V)$. 
Another general class of Banach spaces having property $(V)$ is that of Banach spaces which 
are $M$-ideal in their bidual, {\it i.e.} those for which the canonical decomposition of their 
third dual is an $\ell_1$ decomposition: 
\begin{displaymath}
X^{\ast\ast\ast}= X^\ast \oplus_1 X^\perp
\end{displaymath}
(see \cite{GS1, GS2}). Note that every subspace of a Banach space $M$-ideal of its bidual is 
itself $M$-ideal of its bidual; hence every such subspace has property $(V)$.
\par
On the contrary, a non-reflexive Banach space that does not contain $c_0$ cannot have 
property $(V)$. In particular, $L^1$ does not have this property. Thus, the $L^p$ spaces have 
$(V)$ for $1< p\leq \infty$, whereas $L^1$ does not have it. For the Orlicz spaces, which are, 
in a natural sense, intermediate between $L^1$ and $L^\infty$, D. Leung \cite{Leung} 
proved, when the dual space is weakly sequentially complete, not only that these Orlicz 
spaces have property $(V)$, but that they actually have the local property $(V)$, \emph{i.e.} all their 
ultrapowers have property $(V)$.\par
D. Leung's proof uses non trivial properties of Banach lattices. In this paper, we shall give an elementary proof 
of the (weaker) result that the Orlicz space $L^\Psi$ has property $(V)$, when the complementary function of 
$\Psi$ satifies the $\Delta_2$ condition.
\medskip

\noindent{\bf Acknowledgement.} This work was made during the stay in Lens, in May--June 
2005, of the fourth-named author, as \emph{Professeur invit\'e} of the \emph{Universit\'e 
d'Artois}.\par
We are very grateful to the referee for having simplified the proof of Theorem~\ref{property V}, making it 
shorter and very more elegant and conceptual, by giving us the statement and the proof of 
Proposition~\ref{resultat abstrait}.

\section{The Morse-Transue space}

In this paper, we shall consider Orlicz spaces defined on a probability space 
$(\Omega,\P)$, that we shall assume non purely atomic.\par
By an Orlicz function, we shall understand that $\Psi\colon [0, \infty] \to [0,\infty]$ 
is a non-decreasing convex function such that $\Psi(0)=0$ and $\Psi(\infty)=\infty$. To avoid 
pathologies, we shall assume that we work with an Orlicz function $\Psi$ having the 
following additional properties: $\Psi$ is continuous at $0$, strictly convex (hence 
\emph{strictly} increasing), and such that
\begin{displaymath}
\frac{\Psi(x)}{x}\mathop{\longrightarrow}_{x\to \infty} \infty.
\end{displaymath}
This is essentially to exclude the case of $\Psi(x) = ax$. 
The Orlicz 
space $L^\Psi(\Omega)$ is the space of all (equivalence classes of) measurable
functions $f\colon \Omega\to\C$ for which there is a constant $C > 0$ such
that
\begin{displaymath}
\int_\Omega \Psi\Big(\frac{\vert f(t)\vert} {C}\Big)\,d\P(t) < +\infty
\end{displaymath}
and then $\Vert f\Vert_\Psi$ (the \emph{Luxemburg norm}) is the infinimum of all possible 
constants $C$ such that this integral is $\leq 1$.\par
To every Orlicz function is associated the complementary Orlicz function 
$\Phi=\Psi^\ast\colon [0,\infty] \to [0,\infty]$ defined by:
\begin{displaymath}
\Phi(x)=\sup_{y\geq 0} \big(xy - \Psi(y)\big).
\end{displaymath}
The extra assumptions on $\Psi$ ensure that $\Phi$ is itself strictly convex.\par
\smallskip

Throughout this paper, we shall assume that the \emph{complementary} Orlicz function satisfies 
the $\Delta_2$ condition ($\Phi\in \Delta_2$), \emph{i.e.}, for some constant $K>0$, and 
some $x_0>0$, we have:
\begin{displaymath}
\Phi(2x)\leq K\,\Phi(x),\hskip 1cm \forall x\geq x_0.
\end{displaymath}
This is usually expressed by saying that $\Psi$ satisfies the $\nabla_2$ condition 
($\Psi\in \nabla_2$). This is 
equivalent to say that for some $\beta>1$ and $x_0>0$, one has 
$\Psi(x)\leq \Psi(\beta x)/(2\beta)$ for $x\geq x_0$, and that implies that 
$\frac{\Psi(x)}{x}\mathop{\longrightarrow}\limits_{x\to \infty} \infty$. In particular, 
this excludes the case $L^\Psi = L^1$.\par\smallskip

When $\Phi$ satisfies the $\Delta_2$ condition, $L^\Psi$ is the dual space of $L^\Phi$.
\par\medskip

We shall denote by $M^\Psi$ the closure of $L^\infty$ in $L^\Psi$. Equivalently (see \cite{Rao}, 
page 75), $M^\Psi$ is the space of (classes of) functions such that:
\begin{displaymath}
\int_\Omega \Psi\Big(\frac{\vert f(t)\vert} {C}\Big)\,d\P(t)<+\infty,\hskip 1cm 
\forall C>0.
\end{displaymath}
This space is the \emph{Morse-Transue space} associated to $\Psi$, and  
$(M^\Psi)^\ast = L^\Phi$, isometrically if $L^\Phi$ is provided with the Orlicz norm, and 
isomorphically if it is equipped with the Luxemburg norm (see \cite{Rao}, Chapter IV, 
Theorem 1.7, page 110).\par
We have $M^\Psi = L^\Psi$ if and only if $\Psi$ satisfies the $\Delta_2$ condition, and 
$L^\Psi$ is reflexive if and only if both $\Psi$ and $\Phi$ satisfy the $\Delta_2$ condition. 
When the complementary function $\Phi=\Psi^\ast$ of $\Psi$ satisfies it  
(but $\Psi$ does not satisfy this $\Delta_2$ condition, to exclude the reflexive case), we have 
(see \cite{Rao}, Chapter IV, Proposition 2.8, page 122, and Theorem 2.11, page 123):
\begin{equation}\label{M-ideal}
(L^\Psi)^\ast= (M^\Psi)^\ast \oplus_1 (M^\Psi)^\perp,\tag*{($\ast$)}
\end{equation}
or, equivalently, $(L^\Psi)^\ast= L^\Phi \oplus_1 (M^\Psi)^\perp$, isometrically, with the 
Orlicz norm on $L^\Phi$.\par
For all the matter about Orlicz functions and Orlicz spaces, we refer to \cite{Rao}, or to 
\cite{Kras}.
\par\smallskip

It follows from the preceding equation \ref{M-ideal} that $M^\Psi$ is an $M$-ideal in 
its bidual. Hence $M^\Psi$ and all its subspaces have Pe{\l}czy\'nski's property $(V)$ 
(\cite{GS1, GS2}; see also \cite{HWW}, Chapter III, Theorem 3.4, and the end of this paper). 
This result was shown by D. Werner (\cite{Werner}; see also \cite{HWW}, Chapter III, 
Example 1.4 (d), page 105), by a different way, using the ball intersection property 
(in these references, it is assumed moreover that $\Psi$ does not satisfies the $\Delta_2$ 
condition, but if it satisfies it, the space $L^\Psi$ is reflexive, and so the result is 
obvious).
\smallskip

The proof given in \cite{GS1, GS2} of the fact that Banach spaces which are $M$-ideal in their bidual have 
property $(V)$ uses local reflexivity and the notion of  \emph{pseudo-ball}. We are going to give below a 
slightly different proof, which does not use this last notion, and seems to us more transparent. Let us note that, 
however, a stronger property, namely Pe{\l}czy\'nski's property $(u)$, was shown since then to be satisfied by 
the spaces $M$-ideal of their bidual (see \cite{GL} and, in a more general setting, \cite{GKS}; that follows also  
from \cite{Ro}).\par

\begin{theoreme} {(Godefroy-Saab}, \cite{GS1, GS2})]
Every Banach space which is $M$-ideal in its bidual have property $(V)$.
\end{theoreme}

\noindent{\bf Proof.} Assume that $X^{\ast\ast\ast} = X^\ast \oplus_1 X^\perp$ and let 
$T\colon X \to Y$ be a non weakly compact map. By Gantmacher's Theorem, 
$T^{\ast\ast}\colon X^{\ast\ast} \to Y^{\ast\ast}$ is not weakly compact either. This means 
that $T^{(4)}(X^{(4)}) \not\subseteq Y^{\ast\ast}$. Since 
$X^{(4)}= X^{\ast\ast} \oplus (X^\ast)^\perp$ (canonical decomposition of the third dual of 
$X^\ast$), there exists some $u\in (X^\ast)^\perp$, with $\| u\|=1$ such that 
$T^{(4)}(u)\neq 0$. Now the $M$-ideal property of $X$ gives 
$X^{(4)}= (X^\ast)^\perp \oplus_\infty X^{\perp\perp}$. It follows that 
\begin{displaymath}
\| x + a u\| = \max \{ \| x \| , |a|\}, \hskip 5mm \forall x\in X, \forall a\in \C.
\end{displaymath}
\par

By local reflexivity, we can construct a sequence $(x_n)_{n\geq 1}$ in $X$ equivalent to the 
canonical basis of $c_0$ and such that $\| Tx_n \| \geq \delta >0$ for every $n\geq 1$.\par
For that, let $0<\delta < \| T^{(4)} u\|$, $\eps_n >0$ be such that 
$(1 - \eps_n)\| T^{(4)} u\| > \delta$ and 
$\prod_{n\geq 1} (1+\eps_n)\leq 2$, $\prod_{n\geq 1} (1-\eps_n)\geq 1/2$.\par 
Assume that $x_1,\ldots, x_n$ have been 
constructed in such a way that $\| T x_k\| > \delta$ and 
\begin{align*}
\prod_{k=1}^n (1 -\eps_k) \max \{|a_1|,\ldots, |a_n|\} 
\leq \| a_1 x_1 &+\cdots +a_n x_n \| \\ 
&\leq \prod_{k=1}^n (1 +\eps_k) \max \{|a_1|,\ldots, |a_n|\}
\end{align*}
for every scalars $a_1,\ldots,a_n$.\par
Let $V_n$ be the linear subspace of $X^{(4)}$ generated by $\{u, x_1,\ldots,x_n\}$. By 
Bellenot's version of the principle of local reflexivity (\cite{Bell}, Corollary 7), there 
exists an operator $A_n\colon V_n \to X$ such that 
$\|A_n\|$, $\|A_n^{-1}\|$ are less or equal than $(1+\eps_{n+1})$,  
$A_n$ is the identity on the linear span of $\{x_1,\ldots, x_n\}$ and 
\begin{displaymath}
\big|\, \|T^{(4)}u \| - \|TA_n u\|\,\big|\leq \eps_{n+1}\|T^{(4)}u \|. 
\end{displaymath}
If $x_{n+1}= A_nu$, it is now clear that 
\begin{align*}
\prod_{k=1}^{n+1} (1 - \eps_k) \max \{|a_1|,\ldots, |a_{n+1}|\} 
\leq \| a_1 x_1 &+\cdots +a_{n+1} x_{n+1} \| \\ 
&\leq \prod_{k=1}^{n+1} (1 + \eps_k) \max \{|a_1|,\ldots, |a_{n+1}|\}
\end{align*}
for every scalars $a_1,\ldots,a_{n+1}$ and $\| Tx_{n+1} \| > \delta$.\par
Hence
\begin{displaymath}
\frac{1}{2} \max \{|a_1|,\ldots, |a_n|\} \leq \| a_1 x_1 +\cdots +a_n x_n \| 
\leq 2 \max \{|a_1|,\ldots, |a_n|\}
\end{displaymath}
for every scalars $a_1,\ldots,a_n$. Since $\| Tx_n \| > \delta$, this ends the proof.
\hfill$\square$

\section{Pe{\l}czy\'nski's property $(V)$ for $L^\Psi$.}

As we said, the following result is a particular case of that of D. Leung (\cite{Leung}), but we shall give an 
elementary proof.

\begin{theoreme}\label{property V}(\cite{Leung}) 
Suppose that the conjugate function $\Phi$ of $\Psi$ satisfies the $\Delta_2$ condition. Then, 
the space $L^\Psi$ has Pe{\l}czy\'nski's property $(V)$.
\end{theoreme}
\medskip

As it is well-known (and easy to prove), every dual space with Pe{\l}czy\'nski's property $(V)$ 
is a Grothendieck space: every weak-star convergent sequence in its dual is weakly convergent. 
Hence, we have:
\medskip

\begin{corollaire}\label{Grothendieck}
Suppose that the conjugate function $\Phi$ of $\Psi$ satisfies the $\Delta_2$ condition. Then 
the space $L^\Psi$ is a Grothendieck space.
\end{corollaire}

\noindent{\bf Proof of Theorem \ref{property V}.} We may assume that $L^\Psi$ is a \emph{real} Banach 
space.\par
The proof arises directly from the two following results, since $E= M^\Psi$ is a Banach lattice having 
property $(V)$ and $L^\Psi = (M^\Psi)^{\ast\ast}$.

\begin{lemme}\label{petit lemme}
Suppose that the Orlicz function $\Psi$ does not satisfy the $\Delta_2$ condition. Then for every sequence 
$(g_n)_n$ in the unit ball of $L^\Psi$, there exist a sequence $(f_n)_n$ in $M^\Psi$ and a positive function 
$g \in L^\Psi$ such that $|g_n - f_n| \leq g$.
\end{lemme}

\begin{proposition}\label{resultat abstrait}
Let $E$ be a Banach lattice that has property $(V)$. Suppose that for every sequence $(x_n^{\ast\ast})_n$ in 
$B_{E^{\ast\ast}}$, there are a sequence $(x_n)_n$ in $E$ and a positive $x^{\ast\ast} \in E^{\ast\ast}$ 
such that $| x_n^{\ast\ast} - x_n| \leq x^{\ast\ast}$. Then $E^{\ast\ast}$ has property $(V)$.
\end{proposition}

\noindent{\bf Proof of Lemma \ref{petit lemme}.} Since, by dominated convergence, 
\begin{displaymath}
\lim_{t\to +\infty} \int_\Omega \Psi\big( |g_n|\, \ind_{\{|g_n| > t \}}\big)\,d\P = 0,
\end{displaymath}
we can choose, for every $n\geq 1$, a positive number $t_n$ so big that:
\begin{displaymath}
\int_\Omega \Psi\big(|g_n|\, \ind_{\{|g_n| > t_n \}}\big)\,d\P  \leq \frac{1}{2^n}\raise0,5mm\hbox{,}
\end{displaymath} 
and, moreover such that:
\begin{displaymath}
\sum_{n=1}^{+\infty} \P(|g_n| > t_n) < +\infty.
\end{displaymath} 
This last condition implies, by Borel-Cantelli's lemma, that, almost surely, $|g_n| \leq t_n$ for 
$n$ large enough. Equivalently, by setting:
\begin{displaymath}
 \tilde g_n = g_n\, \ind_{\{|g_n| > t_n \}},
\end{displaymath}
 we have, almost surely $\tilde g_n =0$ for $n$ large enough. It follows that almost surely 
$\sup_n |\tilde g_n|$ is attained. Set now:
\begin{displaymath}
A_n=\{ \omega\in \Omega\,;\ |\tilde g_1(\omega)|,\ldots, |\tilde g_{n-1}(\omega)| 
< |\tilde g_n (\omega)|\ \text{and}\ |\tilde g_k(\omega)|\leq |\tilde g_n (\omega)|,\ \forall k\geq n\}
\end{displaymath}
($\omega\in A_n$ if and only if $n$ is the first time for which $\sup_k |\tilde g_k (\omega)|$ is attained).\par
The sets $A_n$ are disjoint and 
\begin{displaymath}
\sup_{n\geq 1} |\tilde g_n| = \sum_{n=1}^{+\infty} |\tilde g_n|\,\ind_{A_n}.
\end{displaymath}
Hence, if we set:
\begin{displaymath}
g = \sup_{n\geq 1} |\tilde g_n|,
\end{displaymath}
we have $g\in L^\Psi$, since, using the disjointness of the $A_n$'s:
\begin{displaymath}
\int_\Omega \Psi(g)\,d\P = 
\sum_{n=1}^{+\infty} \int_{A_n} \Psi(|\tilde g_n|)\,d\P 
\leq \sum_{n=1}^{+\infty} \int_\Omega \Psi(|\tilde g_n|)\,d\P 
\leq  \sum_{n=1}^{+\infty} \frac{1}{2^n}= 1.
\end{displaymath}
That proves the lemma,  by taking $f_n = g_n - \tilde g_n$, which is in $L^\infty \subseteq M^\Psi$. 
\hfill $\square$
\medskip

\noindent{\bf Proof of Proposition~\ref{resultat abstrait}.} Suppose that $T\colon E^{\ast\ast} \to Y$ is not 
weakly compact. Then there exists a sequence $(x_n^{\ast\ast})_n$ in $B_{E^{\ast\ast}}$ such that 
$(Tx_n^{\ast\ast})_n$ is not relatively weakly compact. Choose $(x_n)_n$ and $x^{\ast\ast}$ as in the 
statement of the Proposition, and set $y_n^{\ast\ast} = x_n^{\ast\ast} - x_n$ for all $n$. We have either:\par
$(a)$ $(Tx_n)_n$ is not weakly compact, or\par
$(b)$  $(Ty_n^{\ast\ast})_n$ is not weakly compact.\par
If $(a)$ holds, $T_{\mid E} \colon E \to Y$ is not weakly compact; hence $T_{\mid E}$ fixes a copy of $c_0$.
\par
If $(b)$ holds, let $I$ be the closed lattice ideal generated by $x^{\ast\ast}$ in $E^{\ast\ast}$, normed so 
that $[- x^{\ast\ast}, x^{\ast\ast}]$ is the unit ball, and let $i \colon I \to E^{\ast\ast}$ be the inclusion map. 
Since $(y_n^{\ast\ast} )_n$ lies in $[- x^{\ast\ast}, x^{\ast\ast}]$, $T \circ i$ is not weakly compact. But 
$I$ is lattice isomorphic to a $C(K)$ space, and hence has property $(V)$. Thus $T \circ i$ fixes a copy of 
$c_0$. So $T$ fixes a copy of $c_0$. \hfill $\square$
\bigskip

\noindent{\bf Remark.} We cannot expect that, for $t_n$ big enough, the functions 
$\tilde g_n$ could have a small norm. For example, let $G$ be a standard gaussian random 
variable ${\cal N}(0,1)$. For  $\Psi= \Psi_2$ ($\Psi_2(x)= {\rm e}^{x^2}-1$), we have, for 
every $t>0$:
\begin{displaymath}
\int_\Omega \Psi_2\Big(\frac{|G|\ind_{\{|G|>t\}}}{\varepsilon}\Big)\, d\P 
= \frac{1}{\sqrt{2\pi}}\int_{|x|>t} 
({\rm e}^{x^2/\varepsilon^2}-1)\, {\rm e}^{-x^2/2}\,dx = +\infty
\end{displaymath}
for every $\varepsilon < \sqrt 2$; that means that 
$\Vert G\ind_{\{|G|>t\}}\Vert_{\Psi_2}\geq \sqrt 2$ 
for every $t>0$ (recall that $\Vert G\Vert_{\Psi_2}= \sqrt{8/3}$: see \cite{LQ}, page 31). 

\section{Concluding remarks and questions}

\noindent{\bf 1.} The full D. Leung's result that $L^\Psi$ have the \emph{local property $(V)$}, \emph{i.e.} 
every ultrapower of  $L^\Psi$ have the property $(V)$ (see \cite{CS}) cannot be obtained straightforwardly 
from our proof. Indeed, since $L^\Psi= (M^\Psi)^{\ast\ast}$ is $1$-complemented in every ultrapower of 
$M^\Psi$, it would suffice to prove that every such ultrapower has property $(V)$; but if  
$\big[(M^\Psi)_{\cal U}\big]^\ast$ contains $(L^\Phi)_{\cal U}$ as a $w^\ast$-dense 
subspace, it is bigger. The ultraprower $(L^\Phi)_{\cal U}$ is not exactly known in general. In the particular 
case of $\Psi= \Psi_2$ ($\Psi_2(x)= {\rm e}^{x^2}-1$), we have (\cite{DK}, Proposition 4.1 and 
Proposition 4.2):
\begin{displaymath}
(L^{\Phi_2})_{\cal U} \cong L^{\Phi_2} (\P_{\cal U}) \oplus L^1(\mu_{\cal U}).
\end{displaymath}

However, since $(L^\Psi)^\ast = (L^\Phi)^{\ast\ast} \cong L^\Phi \oplus_1 L^1(\mu)$, all 
the odd duals of $L^\Psi$ can be written
\begin{displaymath}
(L^\Psi)^{(2n+1)} \cong (L^\Psi)^\ast \oplus_1 L^1(\mu_n).
\end{displaymath}
Hence we get that all the even duals of $L^\Psi$ have the property $(V)$.\par
\smallskip

{\bf 2.} We can define the Hardy-Orlicz spaces $H^\Psi$, in a natural way: it is the subspace 
of $L^\Psi$ consisting of the functions on the unit circle $\T=\partial \D$ which have an 
analytic extension in $\D$; equivalently, it is the subspace of $L^\Psi$ whose negative Fourier 
coefficients vanish. In \cite{Bourgain}, J. Bourgain proved that $H^\infty$ has property 
$(V)$. Does $H^\Psi$ have property $(V)$?\par
Note that the answer cannot follow trivially from our Theorem \ref{property V} since $H^\Psi$ 
is complemented in $L^\Psi$ if and only if $L^\Psi$ is reflexive: indeed, the Riesz 
projection from $L^\Psi$ onto $H^\Psi$ is bounded if and only if $L^\Psi$ is reflexive 
(\cite{Ryan}; see \cite{Rao2}, Chapter VI, Theorem 2.8, page 196), and we have:
\medskip

\begin{proposition}\label{complementation}
Assume that $\Psi\in \nabla_2$. Then the Hardy-Orlicz space $H^\Psi$ is complemented in 
$L^\Psi$ if and only if the Riesz projection is bounded on $L^\Psi$. Hence $H^\Psi$ is 
complemented in $L^\Psi$ if and only if $L^\Psi$ is reflexive. 
\end{proposition}

\noindent{\bf Proof.} Only the necessary condition needs a proof. Assume that there is 
a bounded projection $P$ from $L^\Psi$ onto $H^\Psi$. For every $f\in M^\Psi$, and for every 
$g\in L^\Phi$, the translations $t\mapsto f_t$ and $t\mapsto g_t$ are continuous. Hence we 
can define $\tilde P$ by setting:
\begin{displaymath}
\langle \tilde P f, g \rangle = \int_{\T} \langle P(f_t), g_t \rangle \,dt.
\end{displaymath} 
One has $\|\tilde P f\|_\Psi \leq \| P\|\,\|f\|_\Psi$, so that $\tilde P$ is bounded from 
$M^\Psi$ into $L^\Psi$. On the other hand, it is immediate to see that for every trigonometric 
polynomial $f$, one has, if $e_n(x)= \e^{inx}$:
\begin{displaymath}
\tilde P (f)= \sum_{n\in \Z} \hat f(n) \widehat{P(e_n)}(n)\,e_n.
\end{displaymath}
Since $P$ is a projection, we have $P(e_n)=e_n$ for $n\geq 0$; and since $P$ takes its values 
in $H^\Psi$, we have $\widehat{P(e_n)}(k)=0$ for $k<0$; in particular 
$\widehat{P(e_n)}(n)=0$ for $n<0$.\par
We get therefore:
\begin{displaymath}
\tilde P(f)= \sum_{n\geq 0} \hat f(n) e_n,
\end{displaymath}
that is $\tilde P$ is the restriction to $M^\Psi$ of the Riesz projection. Hence the Riesz 
projection is bounded on $M^\Psi$. By taking its bi-adjoint, we get that it is bounded on 
$L^\Psi$.
\hfill$\square$
\medskip

In Ryan's paper (\cite{Ryan}), it is assumed that $\Psi$ is an $N$-function, 
that is $\lim_{x\to 0} \frac{\Psi(x)}{x}=0$. But we may modify $\Psi$ on $[0,1]$ to get an 
$N$-function $\Psi_1$. Since we work on a probability space $(\Omega, \P)$, the new space $L^{\Psi_1}$ 
is equal, as a vector space, to $L^\Psi$, but with an equivalent norm. Hence Ryan's result remains true 
without this assumption.\par
\smallskip
   
Note that, when the probability space $(\Omega,\P)$ is separable, since we have assumed that 
$\Psi \in \nabla_2$, the reflexivity of $L^\Psi$ is equivalent to its separability 
(see \cite{Rao}, Chapter III, Theorem 5.1, pages 87--88).
\medskip

\noindent{\bf 3.} Property $(V)$ allows us to say that $L^\Psi$ looks like $L^p$, $1 < p \leq \infty$. In some 
sense, it may be seen as close to $L^\infty$ when $\Psi\notin \Delta_2$, since it is not reflexive. However, from 
other points of view, it is closer to $L^p$ with $p<\infty$; on the one hand, it is a bidual space; on the other 
hand, one has:

\begin{proposition}\label{Dunford-Pettis}
If $\Psi\in \nabla_2$, then $L^\Psi$ never has the Dunford-Pettis property.
\end{proposition}

\noindent{\bf Proof.} We are actually going to show that $M^\Psi$ does not have the 
Dunford-Pettis property. That will prove the proposition, since 
$L^\Psi = (M^\Psi)^{\ast\ast}$.\par\noindent
Since $\Psi\in \nabla_2$, there is some $\alpha>1$ and some $c>0$ such that 
$\Psi(x)\geq c x^\alpha$. It follows that $L^\Psi\subseteq L^\alpha$ and the natural 
injection $i\colon L^\Psi \to L^\alpha$ is bounded, and hence weakly compact, since 
$L^\alpha$ is reflexive.\par\noindent
Take now an orthonormal sequence $(r_n)_{n\geq 1}$ in $L^2$ with constant modulus equal to 
$1$ (for example, an independent sequence of random variables taking the values $\pm 1$ each 
with probability $1/2$). One has 
$\int_\Omega r_n f\,d\P\mathop{\longrightarrow}\limits_{n \to +\infty} 0$ for every $f\in L^2$. 
By density, this remains true for every $f\in L^1$, and in particular for every $f\in L^\Phi$, 
since $L^\Phi\subseteq L^1$. Therefore, $(r_n)_{n\geq 1}$ weakly converges to $0$ in $M^\Psi$.  
Since $\|r_n\|_\alpha=1$, $\big(i(r_n)\big)_n$ does not norm-converge to $0$, and hence the 
weakly compact map $i\colon M^\Psi \to L^\alpha$ is not a Dunford-Pettis operator. Therefore 
$M^\Psi$ does not have the Dunford-Pettis property.\hfill$\square$
\medskip

A slightly different way to prove this is to use that for every Banach space $X$ which has the 
Dunford-Pettis property and which does not contain $\ell_1$, its dual $X^\ast$ has the Schur 
property (\cite{Fa, PT}; see also \cite{LQ}, Chapitre 7, Exercice 7.2). But $M^\Psi$ does not 
contain $\ell_1$ (because all its subspaces have property $(V)$; or because its dual $L^\Phi$ 
is separable). Hence $L^\Phi$ would have the Schur property. The same argument as above shows 
that is not the case.\par
\medskip

\noindent{\bf 4.} We have required in this paper that the complementary function $\Phi$ 
satisfies the $\Delta_2$ condition. Hence, in some sense, the space $L^\Psi$ is far from $L^1$. 
We may ask what happens when we are in the other side of the scale, namely when $L^\Psi$ is 
close to $L^1$. But if $\Psi$ satisfies the $\Delta_2$ condition, then $L^\Psi= (M^\Phi)^\ast$ and 
$M^\Phi$, being $M$-ideal of its bidual, has property $(V)$, as said in the Introduction. It follows that 
$L^\Psi$ is weakly sequentially complete (and in fact has property $(V^\ast)$), and if we assume that 
$\Phi\notin \Delta_2$ (so as $L^\Psi$ is not reflexive), then $L^\Psi$ does not have property $(V)$. 
\medskip\goodbreak

\bigskip
\vbox{\noindent{\it 
P. Lef\`evre et D. Li, Universit\'e d'Artois,
Laboratoire de Math\'ematiques de Lens EA 2462, 
F\'ed\'eration CNRS Nord-Pas-de-Calais FR 2956, 
Facult\'e des Sciences Jean Perrin,
Rue Jean Souvraz, S.P.\kern 1mm 18,\par\noindent
62\kern 1mm 307 LENS Cedex,
FRANCE \\ 
pascal.lefevre@euler.univ-artois.fr \hskip 3mm -- \hskip 3mm 
daniel.li@euler.univ-artois.fr
\smallskip

\noindent
H. Queff\'elec,
Universit\'e des Sciences et Techniques de Lille, 
Labo\-ratoire Paul Painlev\'e U.M.R. CNRS 8524, 
U.F.R. de Math\'ematiques,\par\noindent
59\kern 1mm 655 VILLENEUVE D'ASCQ Cedex, 
FRANCE \\ 
queff@math.univ-lille1.fr
\smallskip

\noindent
Luis Rodr{\'\i}guez-Piazza, Universidad de Sevilla, Facultad de 
Matematicas, Dpto de An\'alisis Matem\'atico, Apartado de Correos 1160,\par\noindent 
41\kern 1mm 080 SEVILLA, SPAIN \\ 
piazza@us.es\par}
}

\end{document}